\documentclass[5p,twocolumn,10pt,times]{elsarticle}
\usepackage{amsmath}
\usepackage{bm}
\usepackage{wrapfig}
\usepackage{times}
\usepackage{amsmath}
\usepackage{amssymb}
\usepackage{mathptmx}
\usepackage[stretch=10]{microtype}
\usepackage{comment}
\usepackage[hidelinks]{hyperref}
\usepackage{listings,xcolor}
\usepackage{adjustbox}
\usepackage[framed,numbered]{matlab-prettifier}
\usepackage{xcolor}
\usepackage{paralist}
\usepackage{booktabs}
\usepackage[shortlabels]{enumitem}

\usepackage[caption=false,font=normalsize,labelfont=sf,textfont=sf]{subfig}

\usepackage[color=yellow]{todonotes}
\definecolor{lightGreen}{RGB}{173,223,138}
\definecolor{lightBlue}{RGB}{166,206,227}
\definecolor{lightYellow}{RGB}{255,238,170}

\def \bezier {B{\'e}zier}

\def\COMMENT #1 {}

\begin{document}
\baselineskip11pt
\lstset{language=Matlab}
\begin{frontmatter}

\title{High-accuracy mesh-free quadrature for trimmed parametric surfaces and volumes}

\author[1,2]{David Gunderman\corref{cor1}}
\cortext[cor1]{Corresponding author}
\ead{david.gunderman@colorado.edu}
\author[2]{Kenneth Weiss}
\author[3]{John A. Evans}
\address[1]{Department of Applied Mathematics, University of Colorado Boulder, Boulder, CO 80309, USA}
\address[2]{Lawrence Livermore National Laboratory, 7000 East Avenue, Livermore, CA 94550, USA}
\address[3]{Ann and H.J. Smead Department of Aerospace Engineering Sciences, University of Colorado Boulder, Boulder, CO 80309, USA}

\begin{abstract} 
This work presents a high-accuracy, mesh-free, generalized Stokes theorem-based numerical quadrature scheme for integrating functions over trimmed parametric surfaces and volumes. 
The algorithm relies on two fundamental steps:
\begin{inparaenum}[(1)]
  \item We iteratively reduce the dimensionality of integration using the generalized Stokes theorem to line integrals over trimming curves, and
  \item we employ numerical antidifferentiation in the generalized Stokes theorem using high-order quadrature rules.
\end{inparaenum} 
The scheme achieves exponential convergence up to trimming curve approximation error and has applications to computation of geometric moments, immersogeometric analysis, conservative field transfer between high-order curvilinear meshes, and initialization of multi-material simulations. We compare the quadrature scheme to commonly-used quadrature schemes in the literature and show that our scheme is much more efficient in terms of number of quadrature points used.  We provide an open-source implementation of the scheme in MATLAB as part of QuaHOG, a software package for Quadrature of High-Order Geometries. 
\end{abstract}

\end{frontmatter}

\section{Introduction}
\label{sec:intro}
\begin{figure}
  \centering
  \includegraphics[width=1\linewidth]{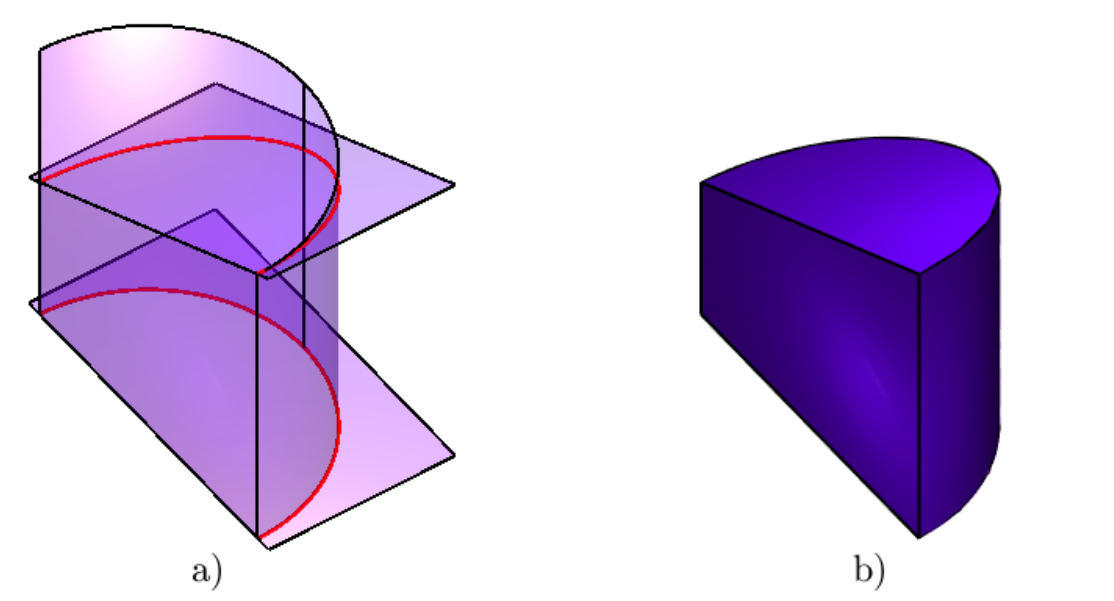}
  \vspace{-.5cm}
  \caption{An example of a region formed by trimmed surface patches. a) The region is bounded by five tensor-product parametric patches 
           with trimming curves shown in red. b) The trimmed object is shown.}
  \label{fig:trimmed_NURBS_example}
\end{figure}

Quadrature over trimmed parametric surfaces and over regions bounded by trimmed parametric surfaces is important in a variety of applications. For example, in computer-aided geometric design (CAGD), three-dimensional geometric objects are typically given in terms of boundary representations (BREPs), often as tensor-product non-uniform rational B-spline (NURBS) surfaces or T-spline surfaces. In many CAGD programs, \emph{trimmed} NURBS patches or T-spline patches-- formed by discarding portions of the patch -- are also allowed 
(c.f. Figure~\ref{fig:trimmed_NURBS_example}) \cite{majeed2017isogeometric,wei2017truncated,wei2018blended,casquero2020seamless,wei2021tuned,wei2021analysis}. 

Some of a trimmed patch's edge curves are formed as intersections of two parametric patches. These edge curves are called \emph{trimming curves} and can be represented either as physical space curves in $\mathbb{R}^3$ or as parametric space curves in the parametric spaces of either trimmed patch~\cite{marussig2018review}. Efficient and accurate integration of material properties over volumes bounded by sets of trimmed patches, or over the surfaces formed as unions of trimmed patches themselves, is often an important operation in the design process, and can be used, for example, to compute an object's volume or surface area~\cite{hafner2019}.

In addition to its application to geometric design, integration over regions bounded by trimmed parametric surfaces is important in the context of computer-aided engineering (CAE), where it is common to use high-order polynomials and/or rational functions as a basis for analysis in the finite element method (FEM). Integration over surfaces and/or volumes formed by trimming typical finite element shapes (such as high-order quadrilateral patches or hexahedral cells) is a fundamental operation in analysis paradigms such as moment-fitting methods \cite{mousavi2010generalized}, immersed finite element methods \cite{burman2010fictitious}, Lagrangian remap methods \cite{anderson2018high}, and multi-mesh methods \cite{johansson2019multimesh}.

For example, in high-order Lagrangian remap methods for computational fluid dynamics (CFD), a curvilinear finite element mesh is allowed to move with the fluid flow until it becomes too distorted for analysis. Then, field variables defined on the distorted source mesh are projected onto an analysis-ready target mesh. Conservation of integrated material quantities, such as total energy, is ensured only with highly accurate integration over the intersections between source mesh elements and target mesh elements \cite{anderson2018high, grandy1999conservative,  margolin2003second}.

\begin{figure*}
  \centering
  \includegraphics[width=1\linewidth]{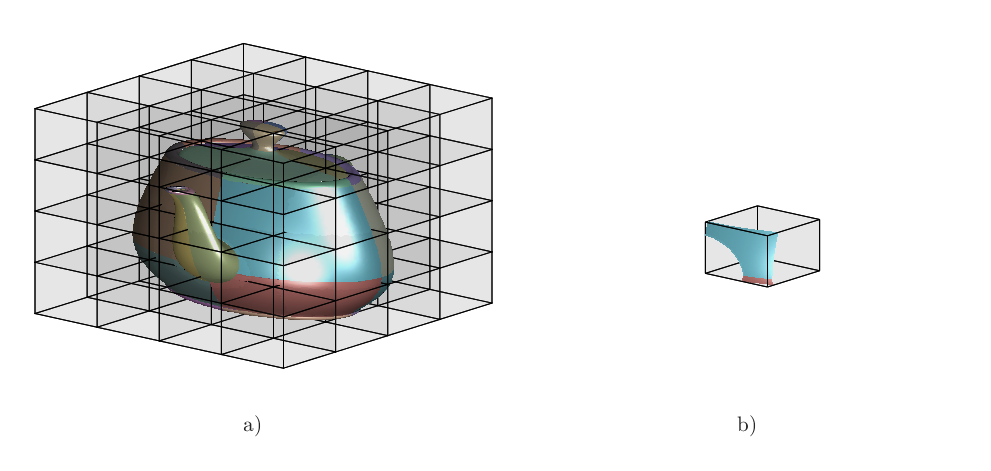}
  \caption{
    \textup{a)} In immersed methods, a shape like the Utah teapot is immersed into a background grid of Cartesian cells on which a finite element space is defined.
    \textup{b)} The boundary is then incorporated by considering where it intersects each cell. In this case, the boundary intersects one of the cells as pictured. 
      Efficient integration must be performed on these types of surfaces, as well as on the volume left after cutting away everything outside the given surface.
  }
  \label{fig:immersed_example}
\end{figure*}

As another example, in immersed finite element methods, an object of geometric interest is immersed in a simpler non-boundary-fitted background mesh. Boundary conditions, domain restrictions, and interface conditions are then taken into account through accurate integration over (1) the volume enclosed by cut cells and (2) the surface of the cut cells \cite{peskin2002immersed, kamensky2015immersogeometric, kudela2015efficient}. See Figure~\ref{fig:immersed_example} for an example of a cut cell formed by a background Cartesian mesh and a geometric object shaped like the Utah teapot.

Efficient high-order integration over typical FEM volume elements and surface elements has been studied extensively for simple cell shapes, such as hexahedral cells and quadrilateral surface patches. Highly optimized rules have been developed in these cases by taking advantage of the existence of a parametric mapping from a simple reference element such as the unit cube or unit square. These methods typically rely on mapping a high order quadrature rule designed for the parametric space into the considered physical space \cite{hughes2010efficient, hughes2012finite}. However, developing high-accuracy quadrature rules for arbitrary cell geometries and topologies, as occur for trimmed surfaces and volumes, is more difficult and is an active area of research.

Accurate integration over arbitrary surfaces and volumes has been approached from multiple standpoints, which can be broadly categorized into (1) domain decomposition methods, such as octree decomposition or body-fitted meshing of the trimmed geometry and (2) dimension-reduction techniques, in which $d$-dimensional integrals are transformed into $d-1$-dimensional integrals. We provide a brief review of these types of methods below. For more complete short surveys of quadrature strategies for arbitrary geometries, see \cite{sudhakar2013quadrature, thiagarajan2014adaptively, divi2020error}.

In domain decomposition-based integration methods, the integration region is decomposed into simpler subregions and then integration is performed over each of the simpler subregions. Domain decomposition can be split into two general categories: low-order domain decomposition and high-order domain decomposition. 
In low-order domain decomposition integration methods, the geometry of interest is decomposed into a set of simple planar-faceted tetrahedra/triangles or hexahedra/quadrilaterals which approximate the original geometry, then standard quadrature rules for these shapes are used to integrate functions~\cite{golub1969calculation, keast1986moderate}.  Because the boundary is approximated by straight-sided geometry, these methods can only achieve low-order convergence to the correct integral with respect to total number of quadrature points. Much work has been done on optimizing these low-order methods~\cite{divi2020error, krishnamurthy2010accurate}. However, the fundamentally low convergence precludes the possibility of high accuracy integration.

Higher-order domain decomposition integration methods rely on a decomposition of the geometric domain into high-order curvilinear elements, such as curved tetrahedra or triangles, and can be thought of as high-order meshing. These methods can attain high-order convergence to the correct integral as the number of subelements increases. However, due to the inherent difficulty of high-order meshing, the pre-processing step of finding a high-order domain decomposition is known to be costly, difficult, and nonrobust \cite{roca2011defining, sherwin2002mesh, engvall2016isogeometric, engvall2017isogeometric, engvall2018mesh, antolin2019isogeometric}. 


\begin{table*}[t]
  \caption{A table of convergence rates under geometric refinement of the approximation of the trimming curves and under integration refinement of the underlying Gaussian quadrature rules using our technique. The refinement strategies are explained further in Section~\ref{sec:quad_ref}.}
  \label{tab:convergence_rates}
  \centering
  \begin{tabular}{ccc}
    \toprule
                            & untrimmed/exact trimming curves & approximating trimming curves\\
    \midrule
       geometric refinement & \emph{N/A} & high-order convergence to integration error\\
     integration refinement & exponential convergence & exponential convergence to geometric error\\
    \bottomrule
  \end{tabular}
\end{table*}

A common method for simplifying quadrature over complicated geometries avoids decomposing the interior by transforming $d$-dimensional integrals over a domain's interior into a sum of $d-1$-dimensional integrals over its boundaries, typically using the generalized Stokes theorem or the implicit function theorem. These methods have been employed extensively for linear polyhedral regions and implicitly-defined regions~\cite{sudhakar2013quadrature,chin2015numerical,saye2015high, olshanskii2016numerical}, particularly when the integrand is symbolically antidifferentiable. To the best of our knowledge, the use of a dimension reduction technique for high-accuracy integration of arbitrary functions over regions bounded by trimmed parametric surfaces or over the surfaces themselves has not been described in the literature, though similar methods have been hinted at for untrimmed surfaces \cite{mousavi2010generalized}.

Additionally, a popular belief about Stokes theorem-based methods is that they can only be used on integrands for which symbolic antiderivatives are available, and these methods have been used in this case for decades~\cite{sheynin2003moment, li1993moment}. However, recently it has become more common to use Stokes theorem-based methods for arbitrary integrands~\cite{sommariva2009gauss,jonsson2017cut} since numerical antiderivatives can be efficiently computed to high precision using a high-accuracy quadrature rule, such as Gaussian quadrature. These methods have, to the best of our knowledge, only recently been used in 2D or surface cases~\cite{gunderman2020spectral}. In this paper, we extend them to the 3D case.

Another interesting approach which has been described for trimmed volumetric domains whose trimming surfaces can be described as implicit functions is that of the error correction term-type schemes, which initially approximates the trimming function with a linear function, then applies a Taylor series correction to obtain higher order convergence \cite{scholz2017first,scholz2019numerical, thiagarajan2014adaptively}. These methods can only be applied to parametric surfaces if they can be implicitized, which can be prohibitively expensive for even quite simple shapes. For example, the intersection of two bicubic surfaces is, in the general case, a 324$^{th}$ degree implicit curve in physical space~\cite{sederberg1984implicit}.

If surfaces are given as polynomial or rational parametric inputs (as is the case for NURBS objects), extra information in the form of a high-order parametrization is available without any extra approximation (e.g., by polyhedra), and this information can be fully exploited to attain the most efficient quadrature rules possible. Similarly, if trimming curves are given as high-order curves, reverting to straight-edged approximations of the trimming curves loses information, which reduces accuracy. 

It is important here to mention the moment-fitting method, a common tool in numerical quadrature. Given accurate moments up to degree $p$ for the geometry, the moment-fitting method produces a quadrature scheme which is \begin{inparaenum}[1)] \item exact for polynomials up to degree $p$ and \item such that quadrature points are at pre-specified locations. \end{inparaenum} Moreover, the quadrature weights of the moment-fitting quadrature rule can be enforced to be positive. Importantly, the moment-fitting method can be combined with any integration procedure and the resulting quadrature rule's accuracy depends on the accuracy of the original integration procedure used for calculation of the moments~\cite{mousavi2010generalized, joulaian2016numerical,hubrich2017numerical}. Our quadrature schemes can be used to accurately calculate the moments which are required as inputs into the moment-fitting algorithm.

In \cite{sommariva2009gauss, gunderman2020spectral}, quadrature schemes with spectral (i.e., faster than any algebraic order) or high-order convergence are developed for area integration over planar domains bounded by polynomial or rational parametric curves, respectively, using an application of Green's theorem combined with numerical antidifferentiation.
In this text, we extend that work to consider two problems. 
  First, we consider the problem of evaluating the surface integral of an arbitrary function over a trimmed parametric surface, $S$. 
  Second, we consider the problem of evaluating the volume integral of an arbitrary function over a region, $V$, bounded by a closed set of trimmed parametric surfaces. 
In both of these cases, we focus on regions that may be difficult to efficiently mesh. Formally, we develop numerical algorithms for computing high-accuracy quadrature schemes for integrals of the forms
\begin{align*}
\iint_{S} f(x,y,z) dS\\
\iiint_V f(x,y,z) dV.
\end{align*}
We assume $V$ is bounded by a set of $n_s$ trimmed tensor product surfaces patches $\{S'_i\}_{i=1}^{n_s}$ that are subsets of untrimmed patches $\{S_i\}_{i=1}^{n_s}$,
\begin{align*}
\partial V &= \cup_{i=1}^{n_s} S'_i:\\
S'_i \subset S_i &= \{x_i(s,t),y_i(s,t),z_i(s,t): 0 \leq s,t \leq 1\}
\end{align*}
We assume $S$ is simply the surface defined by a union of some subset of the trimmed patches $S'_i$. 

The algorithms for producing our proposed quadrature scheme for surface and volume integration rely on two fundamental components:
\begin{enumerate}[(1),topsep=0pt,itemsep=-1ex,partopsep=1ex,parsep=1ex]
  \item We iteratively reduce the dimensionality of integration using the generalized Stokes theorem to line integrals over trimming curves, and
  \item We employ numerical antidifferentiation in the generalized Stokes theorem using exponentially-convergent quadrature rules.
\end{enumerate} 
We note that the presented algorithm can be easily adapted to the case of trimmed and untrimmed simplicial patches, if desired.

The algorithm produces quadrature rules which, for smooth integrands, converge to the correct integral exponentially (i.e., faster than any algebraic order) to machine precision with respect to the number of quadrature points for untrimmed patches and patches with polynomial or rational trimming curves. For trimmed patches with general (e.g. algebraic) trimming curves, the algorithm converges exponentially with respect to the number of quadrature points until the error in the geometric approximation of the trimming curvs dominates the total error.  Under geometric refinement of the trimming curve approximations, the integrals converge with a high order which is determined by the trimming curve approximations and the order of the surfaces. See Table~\ref{tab:convergence_rates} for a summary of these convergence rates. In addition to its convergence properties, our algorithm are computationally efficient. The cost for producing the quadrature rules scales linearly with the number of points in each of the underlying Gaussian quadrature rules,  assuming trimming curves or points along the trimming curve are given \textit{a priori}, as is common in computer-aided design (CAD) systems.

Here we would like to note that the convergence of the strategy is limited if the integrand is not smooth within the bounding box of the domain of integration. If the integrand is $p$-th order continuous, then the algorithm presented will attain $p+1$-th order convergence. Since this strategy is targeted at applications such as immersed boundary methods, where the integrand is smooth within the cut cell, the scheme will converge spectrally.

\begin{figure*}
\begin{center}
\vspace{-.5cm}
\includegraphics[width=.9\linewidth]{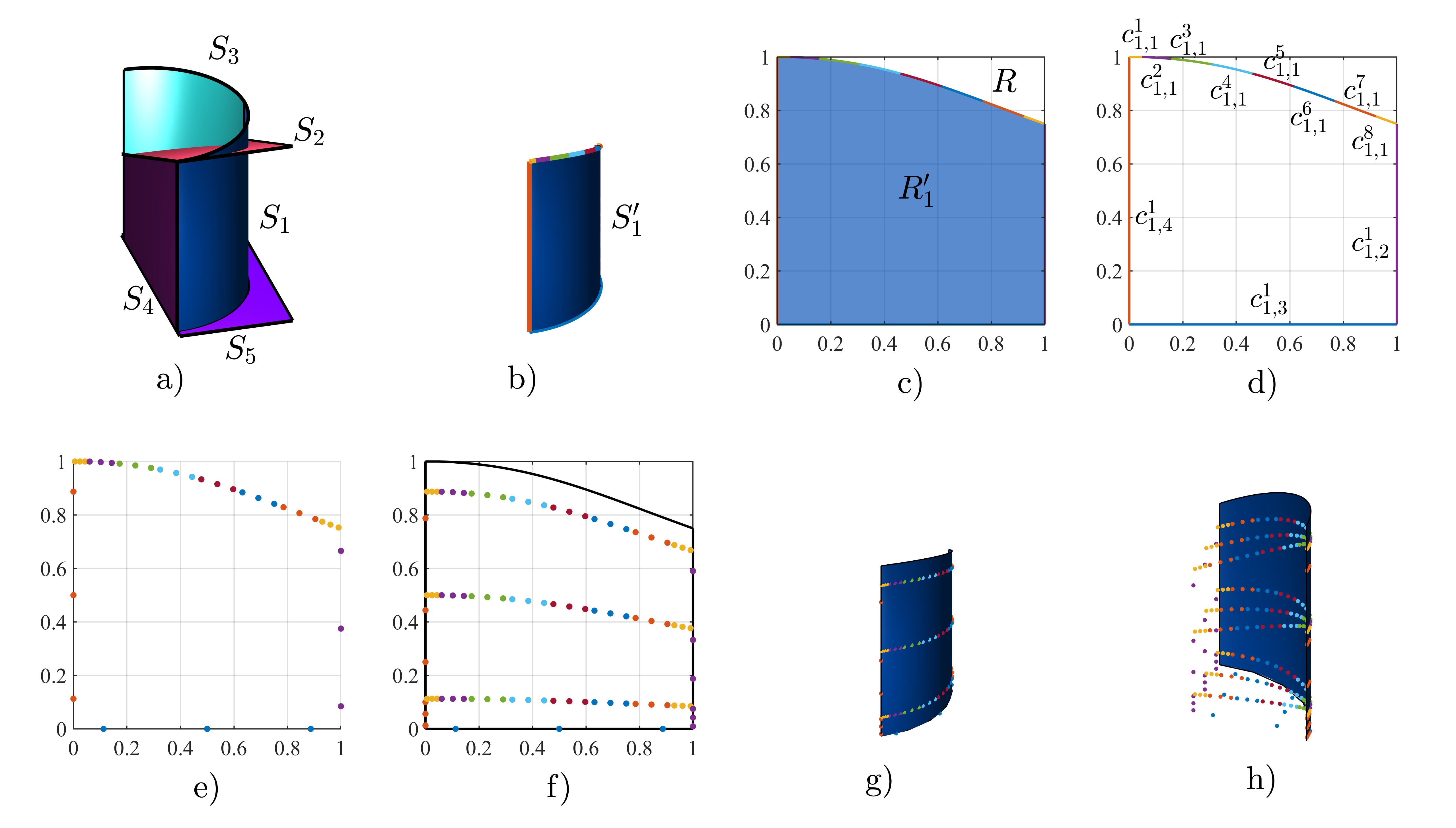}
\vspace{-.5cm}
\caption{A graphical outline of our algorithm. See Section~\ref{sec:algorithm} for an explanation of the notation. a) The integration region is given in terms of its boundaries, which are composed of trimmed surface patches. b) In the case of volume integration, Stokes theorem transforms the integral onto each of the surface patches. In the case of surface integration, the integral is already over the surface patches. c) Each of the surface integrals are then transformed into integrals over the parametric space of the patches, where trimming curves are given (either approximately or exactly) parametrically. d) Green's theorem is applied to transform the parametric area integral onto the trimming curves. e) Integration points are scattered along each trimming curve. f) Numerical Green's theorem is used to scatter quadrature points within the parametric area of interest. g) The parametric quadrature points are mapped back onto the surface. h) Finally, in the case of volume integration, numerical Stokes theorem is used to scatter quadrature points within the original volume.\vspace{-.5cm}
}
\label{fig:graphical_abstract}
\end{center}
\end{figure*}

\paragraph{Contributions}
Among this paper's contributions:
\begin{compactitem}
    \item We present algorithms that compute the locations and weights of high-accuracy quadrature rules for integration of arbitrary smooth integrands over trimmed parametric surfaces and over volumes bounded by such surfaces.
    \item A claimed difficulty of our type of approach in the literature is that it requires the computation of antiderivatives. We show that these are available to high accuracy numerically, easy to calculate for a broad class of useful integrands, and far more efficient than symbolic antiderivatives (see \ref{app:symbolic_numeric}).
    \item We compare our proposed quadrature schemes to three other commonly-used quadrature schemes in the literature that represent the state of the art for integrating over trimmed surfaces or regions bounded by trimmed surface patches. We demonstrate that the proposed algorithms produce quadrature rules which are significantly more efficient in terms of the number of quadrature points needed to achieve a given integration accuracy.
\end{compactitem}
\paragraph{Outline}
The remainder of the paper is organized as follows. 
In Section~\ref{sec:algorithm}, we introduce the new algorithms for producing high-accuracy quadrature rules for trimmed parametric surfaces and regions bounded by such surfaces. 
In Section~\ref{sec:results}, we present numerical results for integrating two integrands over four regions and compare to three other methods used in the literature. In the supplement, we present comparisons for four additional integrands.
 We conclude in Section~\ref{sec:discussion} with a discussion of the results and ideas for future work.

\section{Novel Quadrature Scheme}
\label{sec:algorithm}
In this section, we describe our scheme for producing quadrature rules over trimmed/untrimmed surfaces and volumes. We first show how to compute surface integrals over a set of trimmed surface patches using our scheme.  Then, we show how numerical antidiffentiation can be used to produce quadrature rules for computing volume integrals over a region bounded by a closed union of trimmed/untrimmed surface patches. An example of our algorithm for producing a quadrature rule is shown in Figure~\ref{fig:graphical_abstract}.

\subsection{Surface Integral Case}
\label{sec:surface}
Formally, we consider the problem of finding an appropriate high-accuracy quadrature rule for approximating the surface integral of an arbitrary function $f(x,y,z)$ over a surface $S \in \mathbb{R}^3$ defined as the union of $n_s$ trimmed parametric surface patches, $S = \cup_{i=1}^{n_s} S'_i$:
\begin{equation*}
\iint_S f(x,y,z) dS.
\end{equation*}
 Each trimmed surface patch $S'_i$ is a subset $S'_i \subseteq S_i$ of a corresponding untrimmed tensor-product rational surface patch $S_i$ of bi-degree $(m,n)$:
\begin{align*}
S_i &:= \{(x_i(u,v),y_i(u,v),z_i(u,v)): 0\leq u,v \leq 1\},\text{ where}\\ 
x_i(u,v) &= \frac{\displaystyle\sum_{j=0}^m \sum_{k=0}^{n} w^i_{j,k}x^i_{j,k} B_{j}^m(u) B_{k}^n(v)}{\displaystyle\sum_{j=0}^m \sum_{k=0}^{n} w^i_{j,k} B_{j}^m(u) B_{k}^n(v)} \\
y_i(u,v) &= \frac{\displaystyle\sum_{j=0}^m \sum_{k=0}^{n} w^i_{j,k}y^i_{j,k} B_{j}^m(u) B_{k}^n(v)}{\displaystyle\sum_{j=0}^m \sum_{k=0}^{n} w^i_{j,k} B_{j}^m(u) B_{k}^n(v)}\\
z_i(u,v) &= \frac{\displaystyle\sum_{j=0}^m \sum_{k=0}^{n} w^i_{j,k}z^i_{j,k} B_{j}^m(u) B_{k}^n(v)}{\displaystyle\sum_{j=0}^m \sum_{k=0}^{n} w^i_{j,k} B_{j}^m(u) B_{k}^n(v)},
\end{align*}
where $\{x_{j,k}^i,y_{j,k}^i,z_{j,k}^i\}_{(j=0,k=0)}^{(m,n)}$ are the control points and $\{w_{j,k}^i\}_{(j=0,k=0)}^{(m,n)}$ are the control weights of the $i$-th boundary patch.
In this work, the surface patches are defined in terms of the Bernstein-\bezier\ basis, 
\begin{equation*}
B_{i}^m(u) = \binom{m}{i} (1-u)^{m-i}u^i.
\end{equation*}
Note that another basis or triangular patches could also be used with minimal modifications, but we focus for clarity on the Bernstein-\bezier\, tensor-product patch case. Note also that NURBS surfaces can be easily decomposed into rational tensor-product Bernstein-\bezier\ patches using \bezier-extraction \cite{borden2011isogeometric}, so the algorithm could also be used with minimal modifications for NURBS and related bases (T-splines \cite{bazilevs2010isogeometric}, U-splines \cite{thomas2018u}, hierarchical B-splines \cite{giannelli2012thb}, etc.).

For each component surface patch, $S_i$, there is a set of associated (possibly disconnected) implicitly-defined trimming curves $\{c_{i,k}\}_{k=1}^{m_i} $ in the surface patch's parent parametric space, $R = [0,1]\times [0,1]$, such that $c_{i,k}\subset R$. These trimming curves can be formed by removing the part of $S_i$ on one side of a set of trimming surface patches $\{S_{i_k}\}_{k=1}^{m_i}$, or they may be given \textit{a priori}. Typically included within the set $\{c_{i,k}\}_{k=1}^{m_i}$ are some sections of the original boundary of $R$.  In the parametric space of $S_i$, the set of trimming curves $\{c_{i,k}\}_{k=1}^{m_i}$ splits $R$ into an ``interior" (possibly disconnected) region $R'_i$ and another (possibly disconnected) region $R \setminus R'_i$. The boundary of the interior region $R_i'$ is composed of the trimming curves $\{c_{i,k}\}_{k=1}^{m_i}$. The interior region $R'_i$ maps through the parametric mapping to $S'_i$,
\begin{equation}
S_i' = \{(x_i(u,v),y_i(u,v),z_i(u,v)): (u,v) \in R_i'\}.
\end{equation} See Figure~\ref{fig:graphical_abstract} for an example with notation labeled.

In many cases of interest, the trimming curves $c_{i,k}$ can be represented exactly as a set of parametric curves in the parametric space of $S_i$. In particular, if $c_{i,k}$ is formed by intersecting a planar surface patch $S_i$ and a polynomial (resp. rational) surface patch of bi-degree $(m,n)$, then the trimming curve $c_{i,k}$ will be a polynomial (resp. rational) parametric curve of maximum degree $mn$. In the case that the curves $c_{i,k}$ can not be represented exactly as parametric curves in the space of $S_i$, we assume that the trimming curves $c_{i,k}$ are well-approximated by some set of $n_{i,k}$ pre-computed degree $d_{i,k}$ approximate polynomial or rational parametric curves $\{c^j_{i,k}\}_{j=1}^{n_{i,k}}$ given in the Bernstein-\bezier\ basis and defined in the space $R$ such that $\cup_{j=1}^{n_{i,k}} c^{j}_{i,k} \approx c_{i,k}$:
\begin{align*}
c^j_{i,k} &= \{(u^j_{i,k}(\xi),v^j_{i,k}(\xi)): 0\leq \xi \leq 1 \},\text{ where, }\\ 
u^j_{i,k}(\xi) &= \frac{\displaystyle \sum_{\ell=0}^{d_{i,k}} w^j_{i,k,\ell}u^{j}_{i,k,\ell} B_{\ell}^{d_{i,k}}(\xi) }{\displaystyle \sum_{\ell=0}^{d_{i,k}} w^j_{i,k,\ell} B_{\ell}^{d_{i,k}}(\xi) } \\
v^j_{i,k}(\xi) &= \frac{\displaystyle \sum_{\ell=0}^{d_{i,k}} w^j_{i,k,\ell}v^j_{i,k,\ell} B_{\ell}^{d_{i,k}}(\xi) }{\displaystyle \sum_{\ell=0}^{d_{i,k}} w^j_{i,k,\ell} B_{\ell}^{d_{i,k}}(\xi)} ,
\end{align*}
where $\{(u^{j}_{i,k,\ell},v^{j}_{i,k,\ell}\}_{\ell=0}^{d_{i,k}}$ are the control points and $\{w^j_{i,k,\ell}\}_{\ell=0}^{d_{i,k}}$ are the control weights of the $j$-th approximation curve segment. Our strategy for approximating the trimming curves is given in \ref{app:trimming} and a review of strategies for forming this approximation can be found in \cite{marussig2018review}. If the $c_{i,k}$ curves can be represented exactly using parametric rational curves, then we simply set $n_{i,k}=1$ and use the exact description $c^{1}_{i,k}=c_{i,k}$.  The boundary of the interior region $\partial R_i'$ is then represented, exactly or approximately, by the curves $\{c^j_{i,k}\}_{j=1}^{n_{i,k}}$.


Our algorithm for surface integrals relies on applying numerical Green's theorem in the parametric space $R_i'$ of each trimmed patch to transform each surface integral over $S'_i$ into a sum of line integrals along each approximating boundary curve segment $c^j_{i,k}$.

\subsubsection{Change of Variables to Parent Coordinates}

The algorithm for surface integrals proceeds by considering the integral over each component surface patch individually:
\begin{equation*}
\iint_S f(x,y,z) dS = \sum_{i}^{n_s} \iint_{S_i'} f(x,y,z) dS.
\end{equation*} 
Then, we perform a change of variables to rewrite each of the surface integrals as area integrals in the parametric space of $S_i$:
\begin{equation}
\iint_{S_i'} f(x,y,z) dS = \iint_{R_i'} f(x_i(u,v), y_i(u,v), z_i(u,v)) |\vec{n}_i(u,v)| dudv \label{eq:stokes_parametric}
\end{equation}
where $\vec{n}_i(u,v)$ is the normal vector to $S_i$ at the parametric coordinates $(u,v)$. Note that this can be calculated accurately and efficiently for Bernstein-\bezier\ surface patches using de Casteljau's algorithm \cite{farin2002handbook}. 

\subsubsection{Application of Green's Theorem}
\label{sec:greens}
Once we have moved the problem into the parametric space of $S_i$, we need a highly accurate quadrature rule for integrating an arbitrary function $g(u,v)$ over the region $R_i'$.  See \cite{gunderman2020spectral} for a more detailed derivation of the method we use and a brief survey of other methods. For an arbitrary function $g(u,v)$, by Green's theorem, we have
\begin{equation}
\iint_{R'_i} g(u,v) du dv \approx \sum_{k=1}^{m_i} \sum_{j=1}^{n_{i,k}} \int_{c^j_{i,k}} A_g(u,v) du, \label{eq:greens}
\end{equation}
where we take the $v$ antiderivative,
\begin{equation}
A_g(u,v) = \int_{P_v}^{v} g(u, \xi) d\xi, \label{eq:antiderivative_2d}
\end{equation}
for some appropriately chosen constant $P_v$. Note that the antiderivative can be taken in an arbitrary direction, but we choose the $v$ direction for simplicity of exposition. See Section~\ref{sec:antidiff_const} for a discussion of the choice of antiderivative direction and constant $P_v$. The choice of antiderivative direction is also discussed briefly in \cite{sommariva2009gauss}. From this point forward through the rest of Section~\ref{sec:surface}, we leave off subscript $i,k$ and superscript $j$ for clarity of exposition. The remaining steps are applied to each of the approximating curves which contribute to the sum in Equation~$\eqref{eq:greens}$. 

Since the approximated trimming curves segments are given parametrically, each of the integrals in the sum on the right hand side of Equation~$\eqref{eq:greens}$ can be evaluated in its respective parametric space:
\begin{equation}
\int_{c} A_g(u,v) du = \int_0^1 A_g(u(s),v(s)) \frac{d u(s)}{ds} ds. \label{eq:greens_parametric}
\end{equation}
We again remind the reader that we have suppressed the subscripts and superscripts on the $c$, $u(s)$, and $ v(s)$. The first term in the product can be evaluated with numerical antidifferentiation using a 1D quadrature rule, such as Gaussian quadrature:
\begin{equation*}
A_g(u,v) = \int_{P_v}^{v} g(u, \xi) d\xi \approx \sum_\eta w_\eta g(u,\xi_\eta). \label{eq:gr_antidiff}
\end{equation*} 
%
%
Since we take the parametric mapping $(u(s),v(s))$ to be defined in terms of the Bernstein-\bezier\ basis, the second term in the product can be evaluated stably to machine precision using de Casteljau's algorithm. 
%

\subsubsection{Final Quadrature Rule for Single Surface Patch}

Using a quadrature rule for Equation~$\eqref{eq:greens_parametric}$ yields
\begin{align*}
\int_{c} A_g(u,v) du &= \int_0^1 A_g(u(s),v(s)) \frac{d u(s)}{ds} ds\\
&\approx\sum_{\mu=1}^{m_q} \omega_\mu A_g(u(s_\mu),v(s_\mu)) \frac{d u(s_\mu)}{ds},
\end{align*}
where $\{s_\mu\}_{\mu=1}^{m_q}$ and $\{\omega_\mu\}_{\mu=1}^{m_q}$ are the abcissae and weights of a 1D quadrature rule, respectively, with locations on the interval $[0,1]$. Once these quadrature points and weights have been determined, the antiderivative can also be evaluated using a quadrature rule:
\begin{align*}
A_g(u(s_\mu),v(s_\mu)) &= \int_{P_v}^{v(s_\mu)} g(u(s_\mu), \xi) d\xi\\ & \approx \sum_{\eta=1}^{n_q} w_\eta g(u(s_\mu),\xi_\eta),
\end{align*}
where the $\{\xi_\eta\}_{\eta=1}^{n_q}$ and $\{w_\eta\}_{\eta=1}^{n_q}$ are again the abcissae and weights of a 1D quadrature rule, respectively, with the locations $\xi_\eta$ here spread between $P_v$ and $v(s_\mu)$, where $v(s_\mu)$ will be different for each approximate trimming curve $c_{i,k}^j$ and for each quadrature point $v(s_\mu)$. Note also that the number of quadrature points may be different for each of these antiderivative quadrature rules, but for clarity we use $n_q$ quadrature points for all of them. 
Plugging these quadrature points back in and back-propagating, the approximate form of Equation~$\eqref{eq:greens_parametric}$ becomes
\begin{equation*}
\int_{c} A_g(u,v) du \approx \sum_{\mu=1}^{m_q} \omega_\mu \left(\sum_{\eta=1}^{n_q} w_\eta g(u(s_\mu),\xi_\eta)\right) \frac{d u(s_\mu)}{ds}.
\end{equation*}
The quadrature rule for each component of Green's theorem in Equation~$\eqref{eq:greens}$ (again noting that we have suppressed subscripts and superscripts referring to $i$, $j$, $k$) therefore has formula
\begin{equation*}
\int_{c} A_g(u,v) du \approx \sum_{\mu=1}^{m_q} \sum_{\eta=1}^{n_q} w'_{\mu,\eta} g(u(s_\mu),\xi_\eta),
\end{equation*}
where the quadrature weights are $w'_{\mu,\eta}= \omega_\mu w_\eta \frac{d u(s_\mu)}{ds}$ and the quadrature points are ($u(s_\mu),\xi_\eta)$. In our implementation we use an equal number of quadrature points for each component quadrature rule $n_q=m_q$, but in general they need not be equal.

For each surface patch, we must then multiply the quadrature weights by the normal vector and plug in the quadrature points to the defining parametric equations to convert the area quadrature rule for into a surface quadrature rule,
\begin{align*}
\left(x_{\mu,\eta},y_{\mu,\eta},z_{\mu,\eta}\right) &= \left(x(u(s_\mu),\xi_\eta),y(u(s_\mu),\xi_\eta),z(u(s_\mu),\xi_\eta)\right)\\
w_{\mu,\eta} &= w'_{\mu,\eta}|\vec{n}(u(s_\mu),\xi_\eta)|.
\end{align*}
The quadrature points $\{(u(s_\mu),\xi_\eta)\}$ and weights $\{w'_{\mu,\eta}\}$ correspond to evaluating an area integral corresponding to the right hand side of Eq.~$\eqref{eq:stokes_parametric}$, whereas the quadrature points $\{(x_{\mu,\eta},y_{\mu,\eta},z_{\mu,\eta})\}$ and weights $\{w_{\mu,\eta}\}$ correspond to evaluating surface integrals corresponding to the left hand side of Eq.~$\eqref{eq:stokes_parametric}$.

The quadrature rules for all of the approximated trimming curves are then combined to produce a rule for the trimmed surface patch. Here the explicit dependence on $i$, $k$, and $j$ has been included again for clarity:
\begin{align*}
\iint_{S_i'} &f(x,y,z)dS \\
             &\approx  \sum_{k=1}^{m_i} \sum_{j=1}^{n_{i,k}} \sum_{\mu=1}^{m_q} \sum_{\eta=1}^{n_q} w^j_{i,k,\mu,\eta} f(x^j_{i,k,\mu,\eta},y^j_{i,k,\mu,\eta},z^j_{i,k,\mu,\eta}),
\end{align*}
where $\{(x^j_{i,k,\mu,\eta},y^j_{i,k,\mu,\eta},z^j_{i,k,\mu,\eta})\}$ are the quadrature points and $\{w^j_{i,k,\mu,\eta}\}$ are the quadrature weights.

\subsubsection{Final Quadrature Rule for Full Surface}
Finally, in order to produce a quadrature rule for the entire surface $S$, each of the individual quadrature rules is combined into a full quadrature rule for the original surface:
\begin{align*}
\iint_{S} &f(x,y,z)dS \\
               & \approx \sum_{i=1}^{n_s} \sum_{k=1}^{m_i} \sum_{j=1}^{n_{i,k}} \sum_{\mu=1}^{m_q} \sum_{\eta=1}^{n_q} w^j_{i,k,\mu,\eta} f(x^j_{i,k,\mu,\eta},y^j_{i,k,\mu,\eta},z^j_{i,k,\mu,\eta}),
\end{align*}
where we remind the reader that $n_s$ is the number of surface patches, $m_i$ is the number of trimming curves on patch $S_i$, $n_{i,k}$ is the number of trimming curve approximations for trimming curve $c_{i,k}$, and $m_q$ and $n_q$ are the number of quadrature points used for each trimming curve approximation. Here $\{x_{i,k,\mu,\nu},y_{i,k,\mu,\nu},z_{i,k,\mu,\nu}\}$ are the quadrature points and $\{w_{i,k,\mu,\nu}\}$ are the quadrature weights.
\subsection{Volume Integral Case}
\label{sec:volume_algorithm}
With an accurate quadrature scheme for surface integrals in hand, a quadrature scheme for volume integration can be found by applying an additional iteration of the generalized Stokes theorem and numerical antidifferentiation.

Formally, in this section we consider the problem of finding an appropriate quadrature rule for approximating the volume integral of an arbitrary function $f(x,y,z)$ over a region $V \in \mathbb{R}^3$ bounded by $\partial V = S$, where $S$ is defined as the union of $n_s$ trimmed parametric surface patches, $S = \cup_{i=1}^n S'_i$ as in the previous section. However, in this section we make the additional assumption that $S$ satisfies the Jordan-Brouwer Separation theorem. In other words, the surface $S$ divides $\mathbb{R}^3$ into an inside region $V$ bounded by $S$ and an unbounded outside region $\mathbb{R}^3\setminus V$. We note that in applications, the separation condition is never formally satisfied, since machine precision causes non-water-tightness \cite{marussig2018review}. However, our method is robust to small gaps in models, since we transform integrals onto trimming curves on each surface patch. We wish to compute the integral
\begin{equation*}
\iiint_V f(x,y,z) dV,
\end{equation*}
for arbitrary functions $f(x,y,z)$.
We first use the generalized Stokes theorem to transform this volume integral into a sum of surface integrals. For an arbitrary function $f(x,y,z)$, we have
\begin{align}
\iiint_V f(x,y,z) dV &= \sum_{i=1}^{n_s} \iint_{S_i'} \vec{A_f}(x,y,z) \cdot \vec{dS}, \label{eq:Stokes}
\end{align}
where
\begin{equation}
\vec{A_f}(x,y,z)= \langle 0,0,\int_{P_z}^{z} f(x,y,\xi) d\xi\rangle 
\end{equation}
for some appropriately chosen constant $P_z$.  We discuss the choice of $P_z$ in Section~\ref{sec:antidiff_const}. Note that the antiderivative can be taken in any direction, but we choose the $z$-direction for simplicity of exposition. Because we evaluate the antiderivative in the $z$-direction, the components of $\vec{dS}$ in the $x$- and $y$-directions are zeroed out in the dot product. Equation~$\eqref{eq:Stokes}$ simplifies to the integral of a scalar-valued function,
\begin{equation}
\iiint_V f(x,y,z) dV = \sum_{i=1}^{n_s} \iint_{S_i'} A_f(x,y,z) dS, \label{eq:Stokes_scalar}
\end{equation}
where 
\begin{equation}
A_f(x,y,z)= \int_{P_z}^{z} f(x,y,\xi) d\xi \label{eq:antiderivative_3d} 
\end{equation}
and the normal vector factor $\vec{n_i}$ in the change of variables to parent coordinates in Eq.~$\eqref{eq:stokes_parametric}$ can be replaced with the $z$-coordinate of the normal vector $\vec{n_i}^z$.
Each of the integrals in the sum in Equation~$\eqref{eq:Stokes_scalar}$ can be evaluated to high-precision using a quadrature rule generated with the procedure given in Section~\ref{sec:surface}, with the above caveat:
\begin{equation*}
\iint_{S_i'} A_f(x,y,z) dS \approx \sum_{\sigma=1}^{n_\sigma} w_{i,\sigma} A_f(x_{i,\sigma},y_{i,\sigma},z_{i,\sigma}).
\end{equation*}
Here we have combined the indices $k,j,\mu,\eta$ into one index $\sigma$, where $n_\sigma = m_i n_{i,k} m_q n_q$.
By evaluating the antiderivative numerically using a quadrature rule, we obtain a quadrature rule for each surface patch's contribution to Stokes theorem:
\begin{equation*}
\iint_{S_i'} A_f(x,y,z) dS \approx\sum_{\sigma=1}^{n_\sigma} w_{i,\sigma} \sum_{\psi=1}^{n_p} w_{\psi} f(x_{i,\sigma,\psi},y_{i,\sigma,\psi},z_{i,\sigma,\psi}).
\end{equation*}
The final quadrature rule for the volume case is simply a sum of the contributions from each surface patch:
\begin{equation*}
\iiint_V f(x,y,z) dV \approx \sum_{i=1}^{n_s} \sum_{\eta=1}^{n_\sigma} w_{i,\sigma} \sum_{\psi=1}^{n_p} w_{\psi} f(x_{i,\sigma},y_{i,\sigma},z_{i,\sigma,\psi}),
\end{equation*}
where $\{x_{i,\sigma},y_{i,\sigma},z_{i,\sigma,\psi}\}$ are the quadrature points and $\{w_{i,\sigma}w_\psi\}$ are the quadrature weights in the final quadrature rule.

\subsection{Implementation Details}
Some of the steps of the algorithms above require some choices about implementation which affect the final quadrature scheme, including the choice of antidifferentiation constants. In addition, the quadrature scheme can be simplified if the surface patch of interest is untrimmed. We discuss such implementation details below.
\subsubsection{Choice of Antidifferentiation Constants}
\label{sec:antidiff_const}
In Sections~\ref{sec:volume_algorithm}~and~\ref{sec:greens} of the algorithm given above, the antidifferentiation constants $P_z$ and $P_v$ must be chosen. These antidifferentiation constants should be chosen to be close to the region of integration to ensure that the quadrature points lie close to the domain of integration. In our implementation of the algorithm we choose $P_z$ to be equal to the smallest $z$-value of all of the control points of the surface patches,
\begin{equation*}
P_z = \min_{i,j,k} z^i_{j,k}.
\end{equation*}
This ensures that all of the final quadrature points will lie in the bounding box of the object of interest.

For the constant $P_v$, we choose to use $P_v=0$, since we know \textit{a priori} that all of the trimming curves lie within the original parametric surface patch, whose domain is $[0,1]\times[0,1]$, and this will ensure that all of the surface quadrature points will lie on the original untrimmed parametric surface patch.

\subsubsection{Simplification for Untrimmed Surface Patches}
\label{sec:untrimmed_simplification}
If the surface patch of interest is untrimmed, then the algorithm in Section~\ref{sec:surface} can be simplified by using a tensor-product quadrature rule for the whole surface patch instead of transforming the surface integral into line integrals and using a quadrature rule on each line integral. If there are $n$ quadrature points in the chosen 1D quadrature rule, a tensor-product quadrature has $n^2$ points, whereas a naive implementation of our strategy produces a quadrature with $4n^2$ quadrature points and the same accuracy, including $n^2$ for each of the boundary curves of the surface patch. We note that the left and right curves have geometric correction factors of $\frac{du}{ds} =0$ and the bottom quadrature rule has weights $w_\mu=0$, since the length of the interval in the antiderivative (cf. Eq.~$\eqref{eq:antiderivative_2d}$) is $v(s_\mu)-C_v = 0$. Thus, if we remove the quadrature rules for the left, right, and bottom boundary curves (which contribute nothing to the sum), then our method is equivalent to the tensor-product quadrature rule.  In our implementation, we simply use a tensor-product quadrature rule on the untrimmed surface patches and which is equivalent to using our method with the above simplification.

\section{Results}
\label{sec:results}
In this section, we test the algorithms presented in Section~\ref{sec:algorithm} to calculate surface integrals and volume integrals of some representative functions over various domains defined by trimmed or untrimmed rational or polynomial parametric surface boundaries and their interiors. We first describe the test cases, test integrands, comparison methods, and refinement strategy. Then, we report results for surface and volume integrals. Finally, we consider how our algorithm performs when the trimming curve approximation is refined.
\subsection{Test Domains}
\label{sec:test_domains}
We investigate the method on four datasets: 
\begin{compactenum}[(1)]
  \item A closed propeller shape with no trimming curves, given by 5,136 rational bicubic Bernstein-\bezier\ surface patches (see Figure~\ref{fig:rotor_images}),
  \item A closed bike frame shape given by 2,460 rational bicubic Bernstein-\bezier\ surface patches and 7 planar trimming patches, creating a total of 64 trimming curves (see Figure~\ref{fig:bike_images}),
  \item A non-closed trimmed subset of a sphere shape given as 3 rational biquintic Bernstein-\bezier\ surface patches with 4 trimming surfaces that remove a cookie-cutter pattern out of the shape, creating a total of 6 trimming curves on the original surface (see Figure~\ref{fig:spherestar_images}), and
  \item The closed Utah teapot, a classic example from the computer graphics community composed of 32 polynomial bicubic surface patches and one planar patch, creating a total of 18 trimming curves (see Figure~\ref{fig:teapot_images}).
\end{compactenum}

In the case of the propeller, there are no trimming curves, so for surface integration we use the simplification given in Section~\ref{sec:antidiff_const}, while for volume integrals, the primary addition is the numerical Stokes step of the volume integration algorithm given in Section~\ref{sec:volume_algorithm}. It is important to note that for this text case, the surface quadrature rule for each patch therefore reduces to classical tensor-product Gaussian quadrature, and shows equivalent errors. Since there is no geometric approximation error, we expect exponential convergence under quadrature rule refinement. We discuss our approach to quadrature rule refinement below.

In the case of the bike frame, the trimming curves can be expressed as polynomial parametric curves in the parametric space of the planar patches, and none of the higher-order surface patches are trimmed. Therefore, no trimming curve approximation is necessary. We again expect exponential convergence as there is no geometric approximation error.

In the case of the trimmed sphere shape, the trimming curves are approximated in parameter space using cubic polynomial parametric curves according to the procedure given in~\ref{app:trimming}, since the trimming curves can not be expressed exactly using polynomial or rational parametric curves. Only surface integrals are computed for this example, since it is non-closed.

Finally, in the case of the Utah teapot, the trimming curves are also approximated in parameter space using cubic polynomial parametric trimming curves according to the procedure given in~\ref{app:trimming}, but in this case the shape is closed, so volume integrals are also computed. For the Utah teapot test case, we also investigate the effect of geometric refinement of the trimming curve approximation on the accuracy of the computed integrals.

All experiments and comparisons were implemented in MATLAB and are provided as part of the QuaHOG library~\cite{quahog}.

\subsection{Test Integrands}
\label{sec:test_functions}

For the three closed shapes, we calculate both surface integrals over the boundary and volume integrals over the interior. For the non-closed trimmed sphere, we calculate only surface integrals. 
As test integrands, we used three polynomials of degrees $2$, $4$, and $6$:
\begin{align*}
p_1(x,y,z) &= x^2 - 3xy + \frac{z}{3}\\
p_2(x,y,z) &= \frac{x^3}{3} + x^2 y z - 3y^2z + z + 3\\
p_3(x,y,z) &= y^5 + z^6 - x^2 y z + xz + 2
\end{align*} 
and three representative nonpolynomial functions 
\begin{align*}
f_1(x,y,z) &= \sqrt{x^2 + y^2 + z^2}\\
f_2(x,y,z) &= e^{x + y + z}\\
f_3(x,y,z) &= \frac{x^2 - y + z^2}{x^2 + y^2 + z^2}.
\end{align*} 
These integrands were chosen because they represent classes of functions which might typically be integrated in applications: polynomials of various degrees, algebraic functions with square roots, smooth algebraic functions, and functions with poles. Each of these integrands is smooth within the domains of interest and within the bounding boxes of the domains of interest.  

In the interest of space, we only include comparison plots for $p_3(x,y,z)$ and $f_3(x,y,z)$ in the main text of this paper; results for the other four integrands with comparisons are provided in the paper's supplementary material. 

\subsection{Quadrature Rule Refinement}
\label{sec:quad_ref}
For refinement of our quadrature rules, we use simple, non-adaptive strategies. We use Gaussian quadrature for each of the component 1D quadrature rules and set the number of quadrature points described in Section~\ref{sec:algorithm} equal to each other (i.e., $m_q=n_q=n_p$ for volumes and $m_q = n_q$ for surfaces) and then increase them to refine the quadrature order, thereby increasing accuracy.  We calculate the ``ground truth" value of the integral by increasing the number of quadrature points per trimming curve until the difference between subsequent integration values is smaller than machine precision and, in the case of the trimmed sphere and Utah teapot test cases, by increasing the refinement of the trimming curve approximations until the difference between subsequent integration values with sufficient quadrature points per trimming curve approximation is on the order of machine precision. We perform geometric refinement by increasing the number of cubic trimming curve approximation segments.

\subsection{Comparison Methods} 
\label{sec:comparison}
To show the relative benefit of our methods, we compare with three domain decomposition-based quadrature schemes currently used in the literature. These three strategies represent the state of the art for integrating over trimmed and untrimmed polynomial or rational parametric surfaces and volumes.
\begin{compactdesc}
  \item[\textsc{Octree}.] Our implementation of octree integration adaptively refines a Cartesian grid wherever the boundary intersects the background grid and then uses 3$^{rd}$-order Gaussian quadrature (8 points) on each leaf of the octree. Refinement is performed by increasing the maximum levels of refinement~\cite{kudela2015efficient}. Note that the octree integration approach is only used for volume integrals, since it is a volumetric approach.
  \item[\textsc{Linear mesh}.] We generate a non-adaptive linear simplicial mesh of the geometry using Gmsh~\cite{geuzaine2009} which discretizes the geometry by meshing with a specified maximum element size. We then use 3$^{rd}$-order tensor-product Gaussian quadrature (8 points for tetrahedra and 4 points for triangles) on each element. The various refinement levels represent independently-generated meshes with different maximum element size~\cite{engwirda2014locally}.
  \item[\textsc{Quadratic mesh}.] Similarly to the linear meshing case, we generate a non-adaptive quadratic polynomial simplicial mesh of the geometry using Gmsh~\cite{geuzaine2009} with a specified maximum element size. We then again use 3$^{rd}$-order tensor-product Gaussian quadrature (8 points for tetrahedra and 4 points for triangles) on each element. The various refinement levels represent independently-generated meshes with different maximum element size~\cite{engwirda2014locally}.
 \end{compactdesc}
   We have not optimized any of the implementations of these algorithms nor have we optimized implementations of our own algorithms, but we expect the broad trends, particularly with respect to differences in convergence rates with respect to number of quadrature points, to hold true in any implementation. In our experiments, using higher-order Gaussian quadrature rules for the comparison strategies did not improve the error much, but yielded a much higher number of quadrature points, whereas using lower-order rules increased error. This implies that the geometric approximation itself is the primary source of error. We also note that for the comparison methods, we generate completely new meshes and octrees for each data point, a comparison which favors these methods, since typically a complete remeshing isn't performed to refine the geometry.

\subsection{Results for Surface Integration}
For each of the domains, our method significantly outperforms the two meshing comparison methods in terms of number of quadrature points needed to attain a given error, while avoiding the complicated decomposition of the ambient space necessary in these comparison methods. Note that we do not compare to the Octree strategy, since it is a volumetric approach. The results for surface integration of functions $p_3(x,y,z)$ and $f_3(x,y,z)$ over the four test cases are shown in Figure~\ref{fig:surface_comparisons} and results for the other four test integrands are included in the supplement.

For the propeller model, our method far outperforms the other three methods, because our method retains geometric exactness, allowing for spectral convergence to the correct integral. When our method attains machine precision, the two meshing methods still have errors of at least $10^{\text{-}5}$ or higher, which is a difference in accuracy of more than $10$ orders of magnitude. For the bike model, our method also far outperforms the other three methods. Our method is again able to retain geometric exactness, because the trimming curves in this case are expressed exactly in terms of polynomial parametric curves, which need not be approximated. This again allows our method to attain spectral convergence to the correct integral. In this case, our method at the best point on the plot attains accuracy improvements of $12$ orders of magnitude. Our method also significantly outperforms the comparison methods at lower levels of accuracy for both of these test cases.

The trimmed sphere and Utah teapot models are more nuanced, because they contain trimming curves which must be approximated. Our trimming curve approximation strategy is described in~\ref{app:trimming}. When computing surface integrals over the trimmed sphere model, our method is able to outperform the meshing strategies without requiring any meshing. In this case our method quickly converges to machine precision. For the Utah teapot test case our method still far outperforms the comparison methods, and converges to the geometric approximation error as expected, which is about $10^{-8}$ for this case. 

\subsection{Results for Volume Integration}
Our method also significantly outperforms the two comparison methods in terms of number of quadrature points needed to attain a given error for volume integration.  The results for volume integration of functions $p_3(x,y,z)$ and $f_3(x,y,z)$ over the three closed test cases are shown in Figure~\ref{fig:volume_comparisons} and results for the other four test integrands are again included in the supplement.

As in the surface integration case, our method far outperforms all three comparison methods for the propeller and bike frame models, because no geometric approximation is necessary, allowing for spectral convergence to the correct integral. For the propeller model our method outperforms the comparison methods by nearly $10^{-12}$ at the best point on the plot and is more accurate than the comparison methods with any number of quadrature points.  For the bike frame model our method also far outperforms the other three methods at high levels of accuracy and is only outperformed by  adaptive octree integration at very low accuracies. 
For the Utah teapot test case we observe similar behavior. Our method outperforms the other methods when high accuracy is needed. 

We expect that our method would perform much better at lower accuracies if we implemented it adaptively. 
However, because we have not implemented adaptivity, the adaptive octree strategy is able to outperform it at very low numbers of quadrature points. We defer investigation of an adaptive version of algorithm to another study. Regardless, at high numbers of quadrature points, our method far outperforms the other methods on all of the test cases for high accuracies.

\subsection{Effect of Trimming Curve Approximation}

It is important to note that in the trimmed sphere and Utah teapot test cases, our method achieves spectral convergence only up to the point when geometric approximation error dominates total error. In the case of surface integration over the trimmed sphere case, the chosen trimming curve approximation gives a geometric approximation error less than machine precision, so increasing the integration order reduces the integration error to below machine precision. For surface integration over the Utah teapot test case, the convergence curve flatlines near $10^{\text{-}6}$. For volume integration, it flatlines near $10^{\text{-}9}$. This begs the question, what is the behavior of the error under more refined geometric approximations of the trimming curves? 

After some experimentation, we found that the convergence of the integration error with respect to geometric refinement is very fast. When cubic trimming curve approximations are used (as is the case in this study), the convergence rate is nearly $4$th order for volume integration over the Utah teapot test case with respect to total number of trimming curve segments, as shown in Figure~\ref{fig:teapot_geom}. The order of convergence is, in general, $(d_{i,k}+1)$, where $d_{i,k}$ is the order of the trimming curve approximations. If a fixed number of quadrature points is used for each trimming curve segment under geometric refinement, the convergence rate is also $(r+1)$-th order with respect to the total number of quadrature points for $r-th$ degree polynomial trimming curve approximations. Because this is a volume integral in three-dimensional space, the order of convergence with cubic trimming curves could be considered to be $12$th order with respect to the number of quadrature points per dimension.  This leads us to conclude that if high-precision is needed and the accuracy attained by increasing integration order is not enough, one step of geometric refinement only on the curves (without complicated mesh-refinement) can improve error drastically.

\begin{figure*}
  \includegraphics[width=\linewidth]{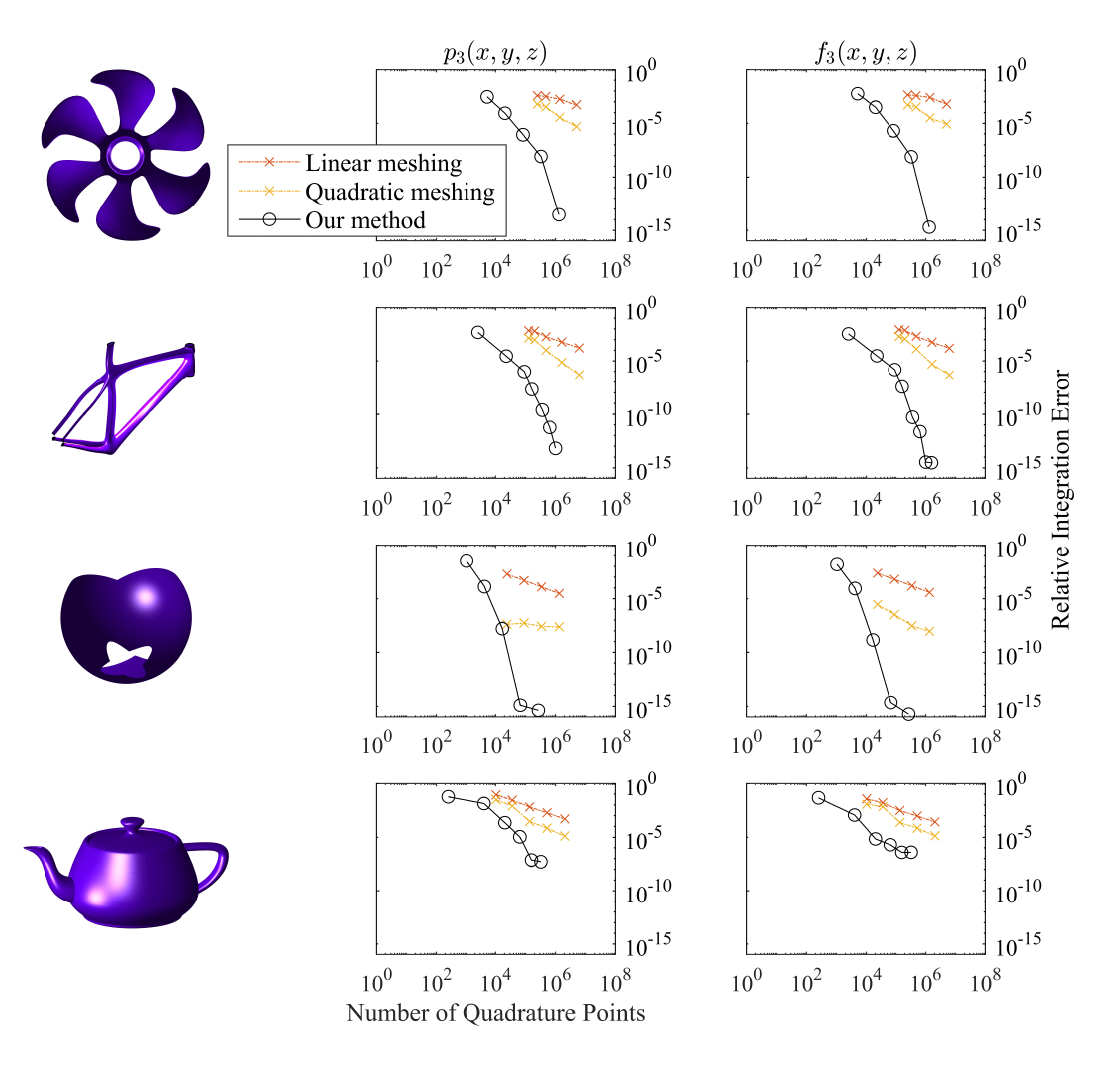}
  \caption{Results of using the algorithm given in Section~\ref{sec:algorithm} to calculate the surface integral of various functions 
          over each of the models defined in Section~\ref{sec:test_domains}. For definitions of the test integrands, see Section~\ref{sec:test_functions}. 
          Note that all plots are on a loglog scale.}
  \label{fig:surface_comparisons}
\end{figure*}

\begin{figure*}
  \includegraphics[width=\linewidth]{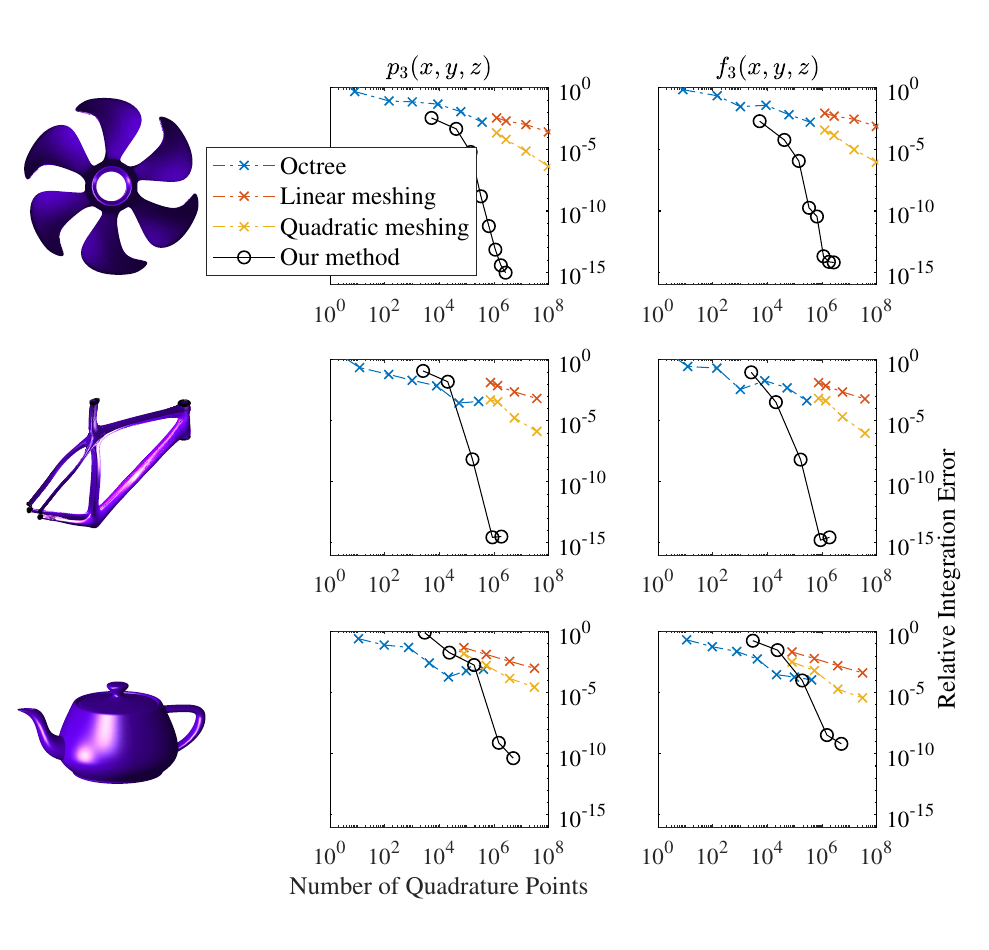}
  \caption{Results of using the algorithm given in Section~\ref{sec:algorithm} to calculate the volume integral of various functions 
          over the interior of each of the closed models defined in Section~\ref{sec:test_domains}. For definitions of the test integrands, see Section~\ref{sec:test_functions}. 
          Note that all plots are on a loglog scale.}
  \label{fig:volume_comparisons}
\end{figure*}

\begin{figure*}
  \centering
  \includegraphics[width=.9\linewidth]{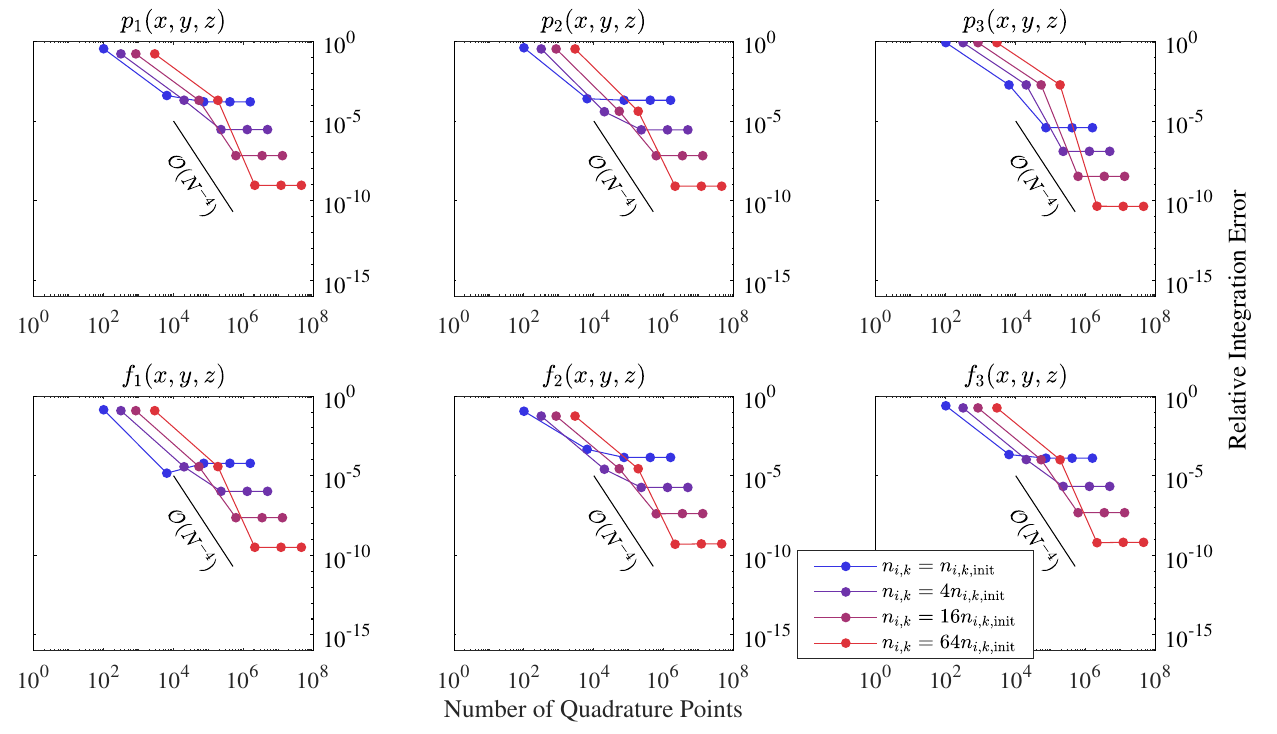}
  \caption{
    Comparing errors due to trimming curve approximation for volume integrals over the interior of the Utah teapot model. 
    Geometric refinement is performed by increasing the number of cubic approximation curve segments $n_{i,k}$ per trimming curve. 
    In each chart, successive data points within a series (colored lines) correspond to increasing quadrature order.
    The set of constants $n_{i,k,\mathrm{init}}$ is the initial number of trimming curve segments used to approximate each of the trimming curves $c_{i,k}$. 
    We note that $n_{i,k,\mathrm{init}}$ in this plot is less than the number of trimming curves used in Figure~\ref{fig:volume_comparisons}. 
    The $N$ inside of the $\mathcal{O}$ notation refers to the total number of quadrature points. 
    For definitions of the test integrands, see Section~\ref{sec:results}. 
    Note that all plots are on a loglog scale.
  }
  \label{fig:teapot_geom}
\end{figure*}

Moreover, for any particular geometric approximation, convergence up to the geometric approximation error is spectral, as shown in Figure~\ref{fig:teapot_geom}. An optimal scheme would switch between using a coarser geometric approximation when lower accuracy is required and finer geometric approximations when higher accuracy is required (i.e., use the bottom-left-most point in each convergence curve of Figure~\ref{fig:teapot_geom}). For the trimmed sphere and Utah teapot test cases in Figures~\ref{fig:surface_comparisons} and \ref{fig:volume_comparisons}, we have simply used a preset geometric approximation and increased the integration order to refine the rule. We leave exploration of the combined geometric and integration order refinement to a future work.

The Utah teapot example shows that the method in question is robust to gaps inherent in B-rep models. Many of the patches are trimmed and these patches do not coincide exactly. However, adjacency information is irrelevant to our algorithm, so the presence of gaps between patches is effectively ignored. These gaps introduce error only in the form of the geometric approximation error already mentioned. Regardless of how large these gaps are, the algorithm itself will function robustly, unlike algorithms based on 3D meshing, octree decomposition, or in/out queries.

\begin{figure*}
 \begin{minipage}{\textwidth}
 \begin{center}
  \includegraphics[width=.35\linewidth]{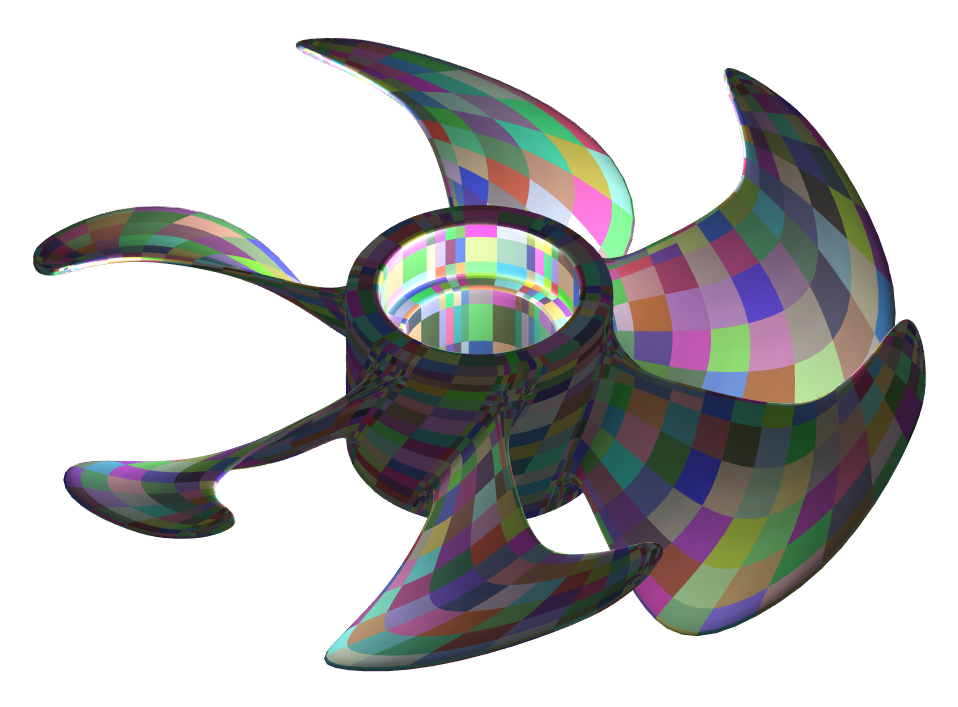}%
  ~\includegraphics[width=.35\linewidth]{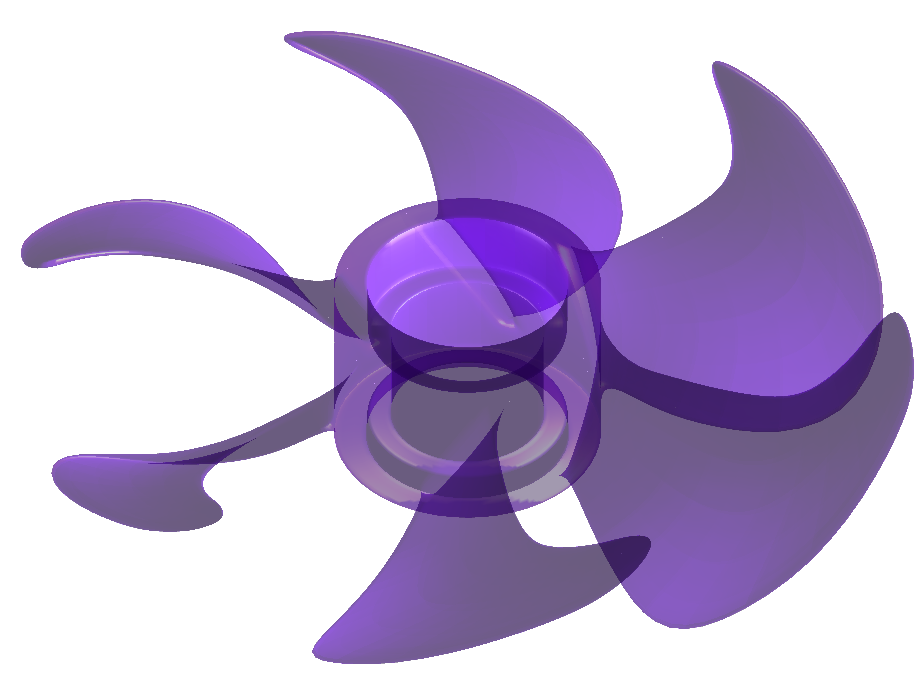}
\end{center}
  \caption{A depiction of the propeller model, which has no trimming curves and 5,136 rational bicubic Bernstein-\bezier\ patches. 
          (Left) An opaque depiction. Colors represent different surface patches of the object. 
          (Right) A partially translucent depiction. 
          }
  \label{fig:rotor_images}
 \end{minipage}\\
 \begin{minipage}{\textwidth}
 \begin{center}
  \includegraphics[width=.35\linewidth]{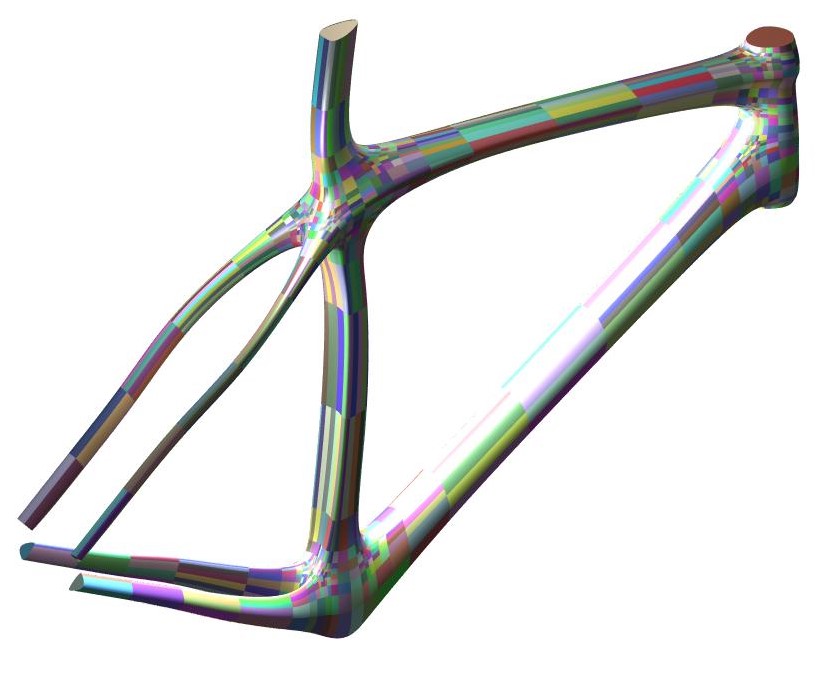}%
  ~\includegraphics[width=.35\linewidth]{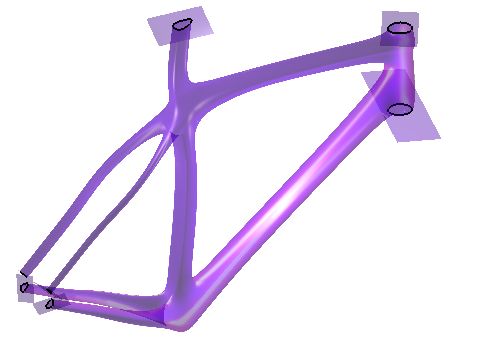}
  \end{center}
  \caption{A depiction of the bike frame model, which has seven trimming curves, including 2,460 rational bicubic Bernstein-\bezier\ patches 
           and seven bilinear Bernstein-\bezier\ planar patches.  
          (Left) An opaque depiction. Colors represent different surface patches of the object. (Right) A partially translucent depiction with trimming curves highlighted in black.}
   \label{fig:bike_images}
 \end{minipage}\\
 \begin{minipage}{\textwidth}
 \begin{center}
  \includegraphics[width=.25\linewidth]{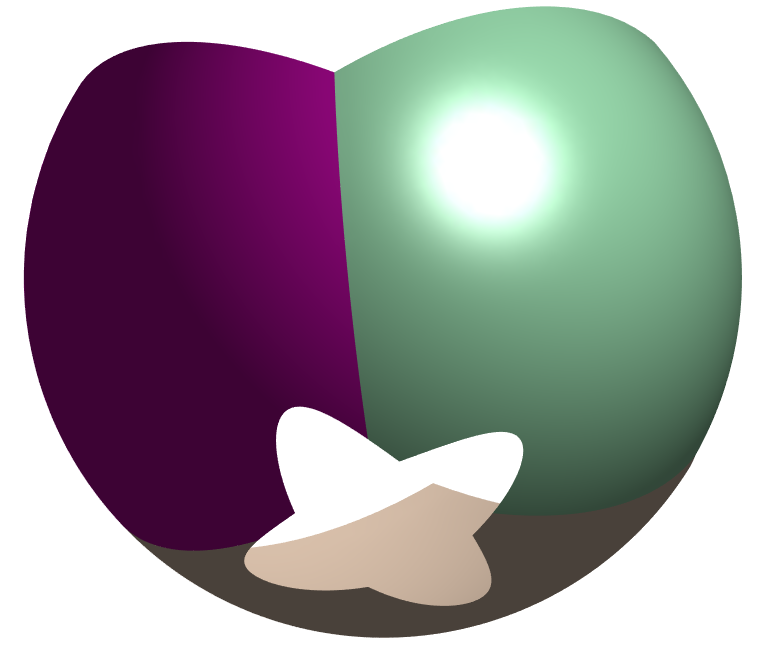}\hspace{.5cm}
  ~\includegraphics[width=.25\linewidth]{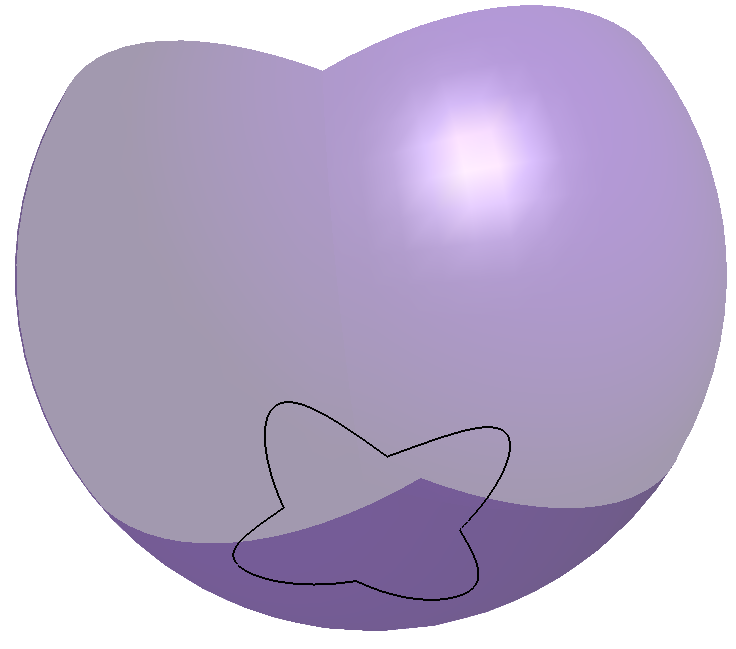}
  \end{center}
  \caption{A depiction of the trimmed sphere model, which has six trimming curves, 
          including three rational biquintic Bernstein-\bezier\ patches. (Left) An opaque depiction. Colors represent different surface patches of the object. (Right) A partially translucent depiction with trimming curves highlighted in black.}
  \label{fig:spherestar_images}
 \end{minipage}\\
 \begin{minipage}{\textwidth}
 \begin{center}
  \includegraphics[width=.31\linewidth]{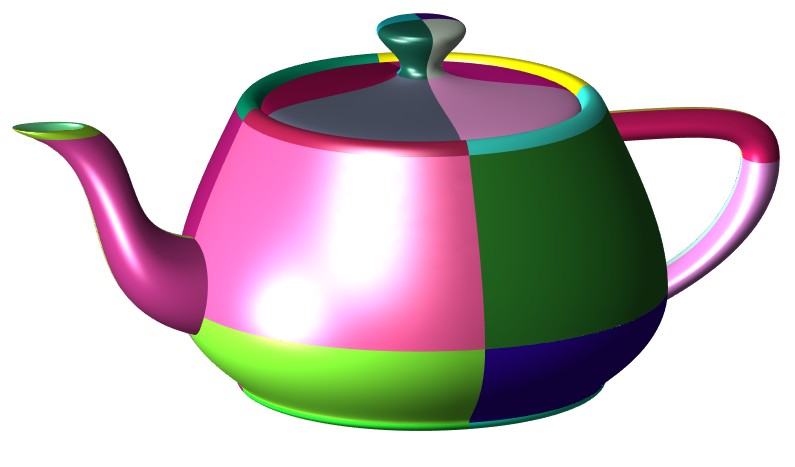}%
  ~\includegraphics[width=.35\linewidth]{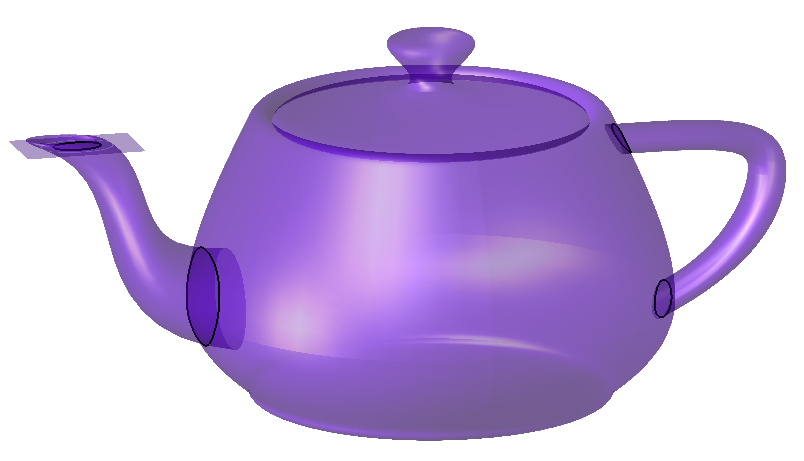}
  \end{center}
  \caption{A depiction of the Utah teapot frame model, which has four trimming curves, 
          including 32 polynomial bicubic Bernstein-\bezier\ patches and one planar patch. 
          (Left) An opaque depiction. Colors represent different surface patches of the object. (Right) A partially translucent depiction with trimming curves highlighted in black.}
  \label{fig:teapot_images}
 \end{minipage}
\end{figure*}

\section{Discussion}
\label{sec:discussion}

In this work, we have presented a high-order, mesh-free quadrature scheme for integrating over trimmed and untrimmed parametric surfaces and volumes. We demonstrated on four shapes and six integrands (three polynomial, three nonpolynomial) that the integration scheme converges spectrally when the surfaces are untrimmed or when the trimming curves can be expressed as low-order polynomial or rational parametric curves. On trimmed surfaces without exact parametric polynomial or rational trimming curves, we showed that the integration scheme attains spectral convergence up to geometric approximation error and high-order convergence when the geometric approximation is refined.

We compared to three different integration schemes: 
\begin{inparaenum}[(1)]
  \item adaptive octree integration, 
  \item linear meshing of the geometry, and
  \item quadratic meshing of the geometry,
\end{inparaenum}
and showed that in comparison to these methods, our method achieves higher accuracy with less quadrature points. We observed the same behavior of our method on simpler shapes that may more commonly appear in applications such as fictitious domain methods, such as the shape in Figure~\ref{fig:trimmed_NURBS_example}.

We would like to briefly note that the high levels of accuracy seen here may seem unnecessary in practice. This may be true, for instance, for certain geometric processing applications such as computation of geometric moments. However, in applications where many integrations must be performed, such as in the solution of time-dependent partial differential equations or in iterative methods, small errors in integration can quickly dominate overall error. High levels of accuracy are also required in applications such as Lagrangian remap, where conservation error is directly proportional to integration error. In addition, we note that while we have focused on high-accuracy integration for our application, the above methods are generally able to achieve lower accuracy with far fewer quadrature points than comparison methods, so these methods could also be applicable in situations where only lower accuracy is necessary. We posit that one of the primary advantages of the methods presented in this paper are that no surface or volume decomposition is necessary for integration, thereby increasing the robustness and decreasing the complexity of the overall integration procedure.

The quadrature schemes presented here can improve accuracy and efficiency in a variety of analysis paradigms, including immersed boundary methods, Lagrangian remap methods, moment-fitting methods, and multi-mesh methods commonly seen in multi-physics applications.  Our plans for future work include three primary thrusts: \begin{inparaenum}[1)]
\item we plan to explore adaptivity of our quadrature scheme by informing the order of quadrature on each trimming curve based on the arc length and curvature of trimming curves, thereby exploring the geometric refinement/integration refinement trade-offs,
\item we plan to investigate spectrally-accurate meshless quadrature schemes which do not have geometric approximation errors, and
\item we also plan to apply these novel high-accuracy quadrature methods to the aforementioned applications.
\end{inparaenum}

\section*{Acknowledgments}
This work was performed under the auspices of the U.S. Department of Energy by Lawrence Livermore National Laboratory under Contract DE-AC52-07NA27344.

\appendix

\section{Trimming curve approximation}
\label{app:trimming}
For surfaces with trimming curves that can not be represented exactly using parametric or rational parametric curves, such as the model in Figure~\ref{fig:trimmed_NURBS_example}, the trimmed sphere example in Figure~\ref{fig:spherestar_images}, and the Utah teapot example in Figure~\ref{fig:teapot_images}, the trimming curves must be approximated. For untrimmed surfaces or surfaces with trimming curves which can be represented exactly, such as the propeller or bike frame models in Figures~\ref{fig:rotor_images}~and~\ref{fig:bike_images}, this step is not required. In many computer-aided design models and engineering applications, the trimming curves are approximated in a pre-processing step before the analysis is started. In these cases, those approximations can be used as the $c_{i,k}^j$ in Section~\ref{sec:surface}. In the case that the trimming curves are not given \textit{a priori}, we must approximate them. 

For each $S_i$, we must approximate each of the $m_i$ trimming curves $\{c_{i,k}\}_{k=1}^{m_i}$ using some set of $n_{i,k}$ polynomial parametric curve segments, $\{c^j_{i,k}\}_{j=1}^{n_{i,k}}$. There is a wide body of literature on producing ordered point sets which lie on trimming curves determined as the intersection of two parametric surfaces and a survey can be found in \cite{marussig2018review}. Once a point set is found, a high-order curve can be fit to the ordered point set, a problem which has also been studied extensively.  In this work, we have used the IRIT polynomial solver \cite{elber2001geometric} to find an ordered, well-spaced point set on an algebraic curve formed as the intersection between two polynomial or rational surfaces. It uses a variety of strategies to achieve this goal, including subdivision strategies and stepping along the curve. Once the ordered point set has been found, we use cubic Lagrange interpolation, with parametric interpolation points spaced based on the physical distance between consecutive points. We also tried using quintic and septic Lagrange interpolation and achieved reasonable results, although the interpolation begins to suffer from the Runge phenomenon as the interpolation order increases, which causes oscillations in high-order polynomial fits to evenly-spaced data.

\section{Symbolic vs Numeric Antidifferentiation}
\label{app:symbolic_numeric}
In the case of nonpolynomial functions, it is often impossible or very difficult to find symbolic antiderivatives of functions, as evidenced for example by the function $f(x,y,z) = \sqrt{x^2+y^2+z^2}$. However, for some classes of functions, such as polynomials, it is relatively easy to find symbolic antiderivatives. In some applications where the integrands of interest are polynomials with known coefficients and are integrated many times, it may be more computationally efficient to simply find antiderivatives symbolically as a pre-processing step. However, numeric antidifferentiation can be preferable even if the integrand is known to be polynomial when the coefficients and/or degree of the polynomial are unknown. This is because symbolically integrating such a polynomial can be cumbersome. In fact, even if the coefficients and degree are known, the symbolic antidifferentiation can still be inefficient because of the relative inefficiency of symbolic operations.

To further explore this dichotomy, we produced $10,000$ polynomials with standard normal distribution of random degrees between $1$ and $20$. We use three strategies to compute the antiderivatives of these functions: 
\begin{inparaenum}[(1)] 
  \item symbolic antidifferentiation using MATLAB's built-in \lstinline{int} function, 
  \item antidifferentiation using a simple numerical implementation of the power rule, and 
  \item numerical antidifferentiation using $10^{th}$ order Gaussian quadrature. 
\end{inparaenum} 
A table with results for pre-processing (i.e., antidifferentiation) of the function along with evaluation (i.e., evaluating at the integration bounds) averaged over $10,000$ random polynomials and $10,000$ random bounds is shown in Table~\ref{tab:sym_num_comp}. As can be seen, symbolic antidifferentiation is far less efficient than a numerical antidifferentiation strategy. The numerical power rule implementation is slightly faster in terms of evaluation time than numerical antidifferentiation using high-order quadrature, but it is far less efficient than Gaussian quadrature in terms of pre-processing time. Importantly, the power rule strategy can only be used in a very specific set of circumstances--on polynomials in the power basis with known coefficients and exponents.
\begin{table}
  \centering
  \caption{A comparison of the average pre-processing and evaluation time for computing the 1D antiderivative of 10,000 polynomials of random degree 
           between $1$ and $20$ with random standard normal coefficients over $10,000$ intervals with random standard normal boundaries. 
           The strategies are described in \ref{app:symbolic_numeric}. Timing is performed in MATLAB.\vspace{.2cm}}
  \label{tab:sym_num_comp}
  \begin{tabular}{cccc}
    \toprule
                            &$(1)$              & $(2)$            & $(3)$              \\
    \midrule
            Pre-processing  & $2\cdot10^{-1}$   & $2\cdot10^{-4}$  & $0$                \\
            Evaluation      & $1.7\cdot10^{0}$ & $8\cdot10^{-7}$  & $2.4\cdot10^{-6}$  \\
    \bottomrule
\end{tabular}
\vspace{-.5cm}
\end{table}

\bibliography{refs}

\end{document}